# REAL POLYNOMIALS WITH A COMPLEX TWIST

**Dr. Bryant Wyatt**
**Department of Mathematics**
**Tarleton State University**
**Box T-0470**
**Stephenville, TX 76402**
**wyatt@tarleton.edu**

**Dr. John Gresham**
**Department of Mathematics**
**Tarleton State University**
**Box T-0470**
**Stephenville, TX 76402**
**jgresham@tarleton.edu**

**Michael Warren**
**Department of Mathematics**
**Tarleton State University**
**Box T-0470**
**Stephenville, TX 76402**
**mwarren@tarleton.edu**

**Abstract**

Student appreciation of a function is enhanced by understanding the graphical representation of that function. From the real graph of a polynomial, students can identify real-valued solutions to polynomial equations that correspond to the symbolic form. However, the real graph does not show the non-real solutions to polynomial equations. Instead of enhancing students' idea of a function, the traditional graph implies a clear disconnect from the symbolic form. In order to fully appreciate the Fundamental Theorem of Algebra, and the non-real solutions of a polynomial equation, traditional graphs are inadequate. Since the early 20th century, mathematicians have tried to find a way to augment the traditional Cartesian graph of a polynomial to show its complex counterpart. Advancements in computer graphics allow us to easily illustrate a more complete graph of polynomial functions that is still accessible to students of many different levels. The authors will demonstrate a method using modern 3D graphical tools such as GeoGebra to create dynamic visualizations of these more complete polynomial functions.

**Introduction**

A good understanding of functions consists of fluency with each representation as well as with the connections between those representations. As students begin to work with the idea of a function, their comprehension can be helped or hindered by the instructor's use of representations (Star & Rittle-Johnson, 2009). In particular, the connections between the symbolic and graphical representation play an important role in students' comprehension of functions (Lipp, 1994; Presmeg, 2014). When the graphical representation of a function demonstrates visually what is being studied symbolically, student comprehension is significantly strengthened. Star and Rittle-Johnson (2009) give an example of a linear relationship between the amount of time spent on a cell phone and the subsequent charges. They concluded that students who use both symbolic and graphical representations to describe the charges had a deeper comprehension of linear functions than students who use only the symbolic form. It is important for students to experience these multiple representations if they are to become robust problem solvers.

Martinez-Planell & Gaisman examined students' understanding of functions of two variables in multivariate calculus (2012). They investigated students' understanding of the graphical representations of functions of two variables. The researchers conclude that students will be better prepared for studying multi-variate functions if they have more chances to develop domains, ranges, and graphical representations of a variety of functions as well as having a good schema for $R^3$. Students need to be proficient with these concepts in order to be well prepared for advanced mathematical study.

One opportunity for instructors to highlight the connections between representations of functions and to confront students with the development of a sophisticated domain, range, and graphical representation of a function is polynomials. Traditionally when visualizing a function with the graphical representation, the domain is selected and those domain values are used as the input values distributed along the horizontal axis. We then take these inputs and run them through the function to produce corresponding outputs that we distribute along the vertical axis. Extending lines vertically and horizontally from these input output pairs, respectively, we place a point at their intersections. We connect these points with line segments to produce a continuous graph. With the proper selection of input values this graph can be a nice two-dimensional approximation to the true function which can be readily produced with pencil and paper or common graphing utilities. This technique extends to functions that map a two dimensional set into a one dimensional set by placing the selected input values on a horizontal plane and the output values they create on a vertical axis perpendicular to this plane. Points are now placed in a three dimensional space where vertical lines extended from the input values intersect horizontal planes extended from the corresponding output value. This creates a three dimensional object that can be depicted with the use of pencil and paper, modern computer graphing programs, or even 3D printed. This technique of visualizing the function breaks down when the sum of the dimensions is greater than three. Visualizing a function that maps complex valued inputs to complex valued outputs is difficult because the resulting object is four dimensional. One way to use the traditional method of graphing such functions is to restrict the domain so that the outputs generated are one dimensional. This creates a three dimensional object which can be visualized. Though not the complete object under study, it is a view which can provide great insight into the function. The idea is not new and was studied in the early 20[th] century but generating three-dimensional graphs by hand was difficult and time consuming, so it was not taught in mainstream mathematics (see Curtis, 1938; Gehman, 1941; Yanosik, 1936; & Yanosik, 1943). Harding and Engelbracht, in a historical look at complex roots, present this idea as an important area of further study

(2007). In an article by Alan Lipp (2001) modern three-dimensional graphing was applied to this topic but is still not widely used. We present the idea here using GeoGebra, a free computer algebra system with dynamic 3D graphing functionality. Our hope is that through interactive graphics, this technique of visualizing polynomials will be more commonly used to help students gain a deeper understanding of polynomial functions.

**The Quadratic**

To begin, consider the simple quadratic $x^2 - 4$. A commonly asked question about this function is, what input or inputs generate an output of zero? Students can verify the answers to this question by looking at the graph, see *figure 1*. The symbolic and commonly used graphical representations are aligned and strengthen a student's understanding of the function. Now consider another simple quadratic $x^2 + 4$. When addressing the same question as before, the symbolic and commonly used graphical representations become disjoint causing a student confusion when attempting to verify the solutions, see *figure 2*. This disconnection between the representations weakens a student's understanding of the function. This problem also occurs for any desired output from this function that is less than four. Hence, traditional graphing techniques are not sufficient and causes confusion for students.

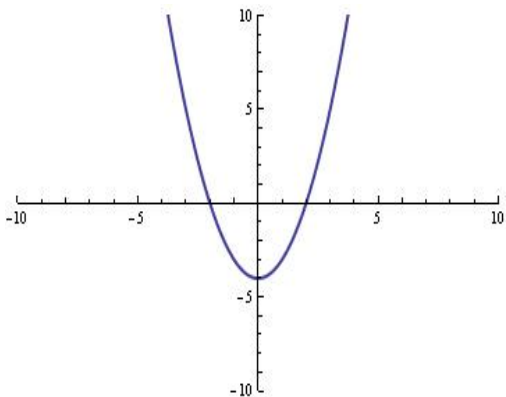
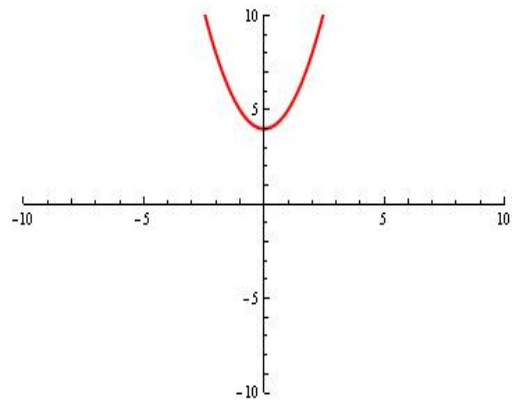

*Figure 1*. Graph of $x^2 - 4$. The two zeros $\pm 2$ are easily verified.

*Figure 2*. Graph of $x^2 + 4$. The two zeros $\pm 2i$ are not easily verified.

To alleviate the confusion with the quadratic $x^2 + 4$, ask what inputs on the complex plane produce an output of zero. However, don't stop with zero. Ask the same question for outputs of 3, 2, 1, -1, etc… and combine these with the traditional outputs greater than or equal to 4. Plot these points above the complex plane and an interesting pattern will emerge. The traditional parabola will be opening up along the real axis but another parabola opening down will appear along the imaginary axis. These parabolas will be connected at their vertices, see *figure 3*. Now all the answers to those commonly asked questions of the quadratic $x^2 + 4$ can be verified. This will also confirm what the students know about the Fundamental Theorem of Algebra. Slice this modified graph at any point with a horizontal plane and you will have two points of intersection except at the vertices which will be a root of multiplicity two. The extension of this idea to an arbitrary quadratic $ax^2 + bx + c$ shifts the vertex of intersection accordingly, see *figure 4*.

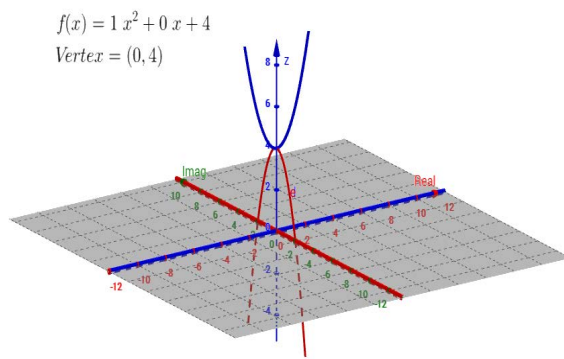
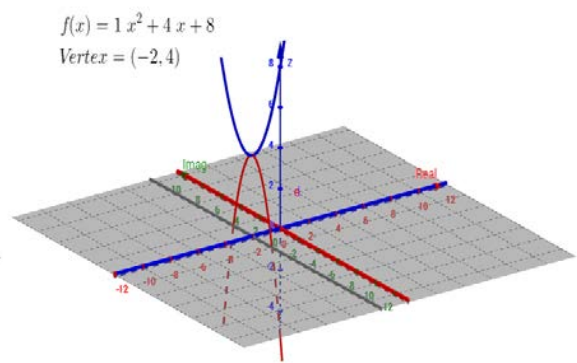

*Figure 3*. 3D Graph of $z^2 + 4$. The two zeros $\pm 2i$ are now easily verified.

*Figure 4*. 3D Graph of $z^2 + 4z + 8$. The non-real parabola lies above the line $x = -\frac{b}{2a} = -2$.

This process restricts the domain of a given function to only those complex valued inputs that produce real valued outputs. The restricted domain of the simple quadratic $x^2 + 4$ can be derived as follows.

Let $f(z) = z^2 + 4$,
where $z = x + iy$ for $x, y \in \mathbb{R}$ and $i^2 = -1$.
Then, $f(z) = (x + iy)^2 + 4$
$= x^2 + i2xy - y^2 + 4$
$= (x^2 - y^2 + 4) + i(2xy)$.

In order to produce a real output, $2xy = 0$.

Therefore,
(1) $y = 0$
or
(2) $x = 0$.

Hence the restricted domain for this quadratic is the real and imaginary axes. Parametric equations in Geogebra can be used to view this graph in three-dimensions (see Appendix). The result is shown in *figure 3*.

Now consider the general quadratic with real valued coefficients.

Let $f(z) = az^2 + bz + c$, ($a, b, c \in \mathbb{R}$ and $a \neq 0$)
where $z = x + iy$ for $x, y \in \mathbb{R}$ and $i^2 = -1$.

Then, $f(z) = a(x + iy)^2 + b(x + iy) + c$

$$= ax^2 + i2axy - ay^2 + bx + iby + c$$
$$= (ax^2 - ay^2 + bx + c) + i(2axy + by).$$

In order to produce a real output,
$$2axy + by = 0$$
$$y(2ax + b) = 0.$$

Therefore,
  (3) $y = 0$
     or
  (4) $x = -\dfrac{b}{2a}.$

Hence the restricted domain for the general quadratic is the real axis and the line parallel to the imaginary axis that lies under the joining vertex of the two parabolas. These perpendicular lines are always centered at the vertex, $= -\dfrac{b}{2a}$ which is where the first derivative of the quadratic is equal to zero.

**The Cubic**
Now apply this process to the general cubic.

Let $f(z) = az^3 + bz^2 + cz + d$, $(a, b, c, d \in \mathbb{R} \text{ and } a \neq 0)$
where $z = x + iy$ for $x, y \in \mathbb{R}$ and $i^2 = -1$.

Then, $f(z) = a(x + iy)^3 + b(x + iy)^2 + c(x + iy) + d$
$$= ax^3 + i3ax^2y - 3axy^2 - iay^3 + bx^2 + i2bxy - by^2 + cx + icy + d$$
$$= (ax^3 - 3axy^2 + bx^2 - by^2 + cx + d) + i(3ax^2y - ay^3 + 2bxy + cy).$$

In order to produce a real output the imaginary part must equal zero,
$$3ax^2y - ay^3 + 2bxy + cy = 0$$
$$y(3ax^2 - ay^3 + 2bx + c) = 0.$$

Therefore,
  (5) $y = 0$
     or
  (6) $3ax^2 - ay^3 + 2bx + c = 0.$

Solving for $y$,
  (7) $y = \pm\sqrt{\dfrac{3ax^2 + 2bx + c}{a}}.$

The restricted domain includes the real axis once again, but (6) seems more complicated. With some work, it can be shown that (6) is a hyperbola in the complex plane.

(8) $\dfrac{\left(x+\frac{b}{3a}\right)^2}{\frac{1}{3}\left(\frac{b^2-3ac}{9a^3}\right)} - \dfrac{y^2}{\left(\frac{b^2-3ac}{9a^3}\right)} = 1$

All cubics have one inflection point, where the second derivative is equal to zero. This hyperbola is centered at the point of inflection, where $x = -\dfrac{b}{3a}$. Using the restricted domain we can now build a 3D graph of any cubic function. The slope at the inflection point will divide all cubics into three distinct categories, see *figure 6*. These categories also correspond to the three possible states of the restricted domain, see *figure 7*. This graph of the cubic also confirms what students know about the Fundamental Theorem of Algebra. Horizontal planes will always intersect this 3D cubic three times except at points where the real and non-real pieces of the cubic touch which will be roots of multiplicity greater than one.

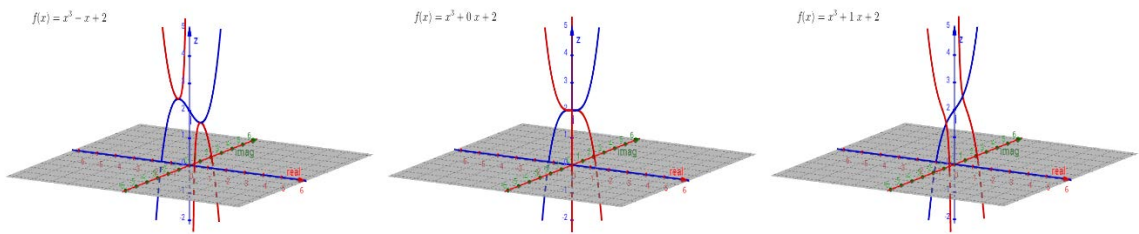

*Figure 6.* 3D Graphs of Cubics. This figure illustrates the three distinct categories of cubic polynomials. From left to right, the slope at the inflection point is negative, zero, and positive.

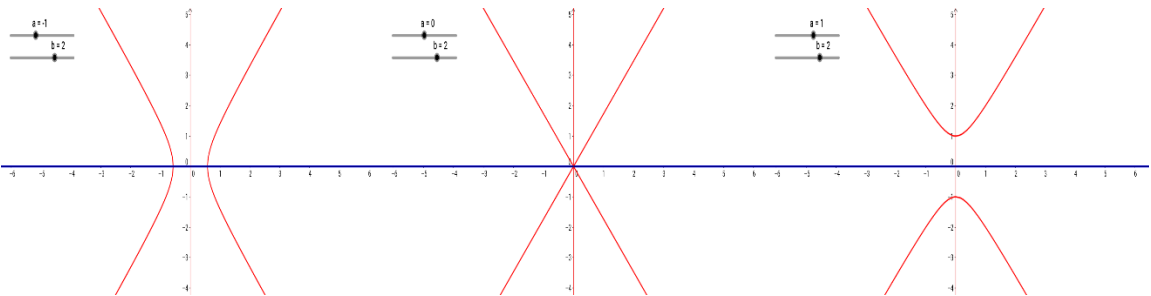

*Figure 7.* Cubic Restricted Domain. This figure illustrates the three distinct categories of the restricted domain for cubic polynomials. From left to right, the slope at the inflection point is negative, zero, and positive.

The quartic acts in much the same way as the cubic but this time its central point will be where the third derivative is equal to zero at the value $x = -\dfrac{b}{4a}$. Let's find its restricted domain. This will take some work.

Let $f(z) = az^4 + bz^3 + cz^2 + dz + e$, $(a, b, c, d, e \in \mathbb{R}$ and $a \neq 0)$
where $z = x + iy$ for $x, y \in \mathbb{R}$ and $i^2 = -1$.

Then, $f(z) = a(x + iy)^4 + b(x + iy)^3 + c(x + iy)^2 + d(x + iy) + e$
$= ax^4 + i4ax^3y - 6ax^2y^2 - i4axy^3 + ay^4 + bx^3 + i3bx^2y - 3bxy^2 - iby^3$
$\quad + cx^2 + i2cxy - cy^2 + dx + idy + e$
$= (ax^4 - 6ax^2y^2 + ay^4 + bx^3 - 3bxy^2 + cx^2 - cy^2 + dx + e) +$
$i(4ax^3y - 4axy^3 + 3bx^2y - by^3 + 2cxy + dy)$.

In order to produce a real output the imaginary part must equal zero,
$4ax^3y - 4axy^3 + 3bx^2y - by^3 + 2cxy + dy = 0$
$y(4ax^3 - 4axy^2 + 3bx^2 - by^2 + 2cx + d) = 0$.

Therefore,
(9)  $y = 0$
or
(10)    $4ax^3 + 3bx^2 + 2cx + d - 4axy^2 - by^2 = 0$.

Solving for $y$,

(11)    $y = \pm\sqrt{\dfrac{4ax^3 + 3bx^2 + 2cx + d}{4ax + b}}$.

The restricted domain for the quartic produces a hyper-hyperbola along with the real axis. One example of the restricted domain and 3D quartic is shown in *figure 8*.

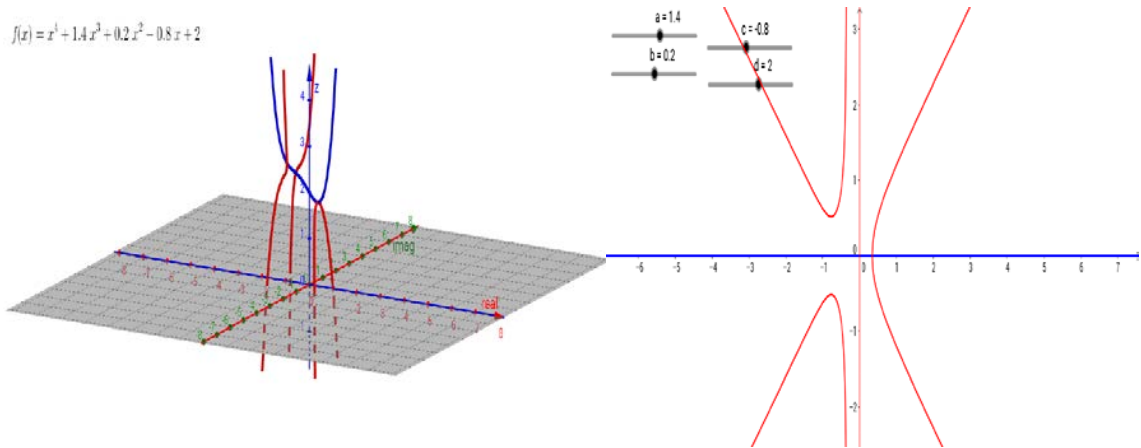

*Figure 8.* This figure illustrates one possible three-dimensional quartic polynomial along with its restricted domain.

The process can continue in this manner to polynomials of higher degree (see note in Appendix). The complexity will increase but the general ideas will remain the same. This technique is not restricted to polynomial functions. It is left to the reader to investigate further.

**Conclusion**

Students need to see the connections between the symbolic and graphical representations of a function. Traditional graphs of polynomial functions are disconnected from what students are told is true with the Fundamental Theorem of Algebra. This creates a confusing paradox that is not necessary with the technology available today. Domain restriction coupled with interactive 3D graphics in programs like GeoGebra alleviates this confusion. Students now have the tools to investigate a more complete polynomial function with symbolic and graphical representations that support each other.

## Appendix

The following is a method to plot the 3D quadratic functions in GeoGebra.

### $x^2 + 4$

From (1), let $x(t) = t$, $y(t) = 0$, and $z(t) = t^2 + 4$.

In GeoGebra input Curve[t, 0, t² + 4, t, ParameterStart, ParameterEnd] where ParameterStart and ParameterEnd give the desired range of values for t.

From (2), let $x(t) = 0$, $y(t) = t$, and $z(t) = -t^2 + 4$.

In GeoGebra input Curve[0, t, -t² + 4, t, ParameterStart, ParameterEnd] where ParameterStart and ParameterEnd give the desired range of values for t.

### $ax^2 + bx + c$

Create sliders for each coefficient value, $a, b$ and $c$. In GeoGebra input Slider[Min, Max, Increment] where Min and Max give the desired range for each coefficient value and Increment gives the quantity to be scrolled with each click on the slider.

From (3), let $x(t) = t$, $y(t) = 0$, and $z(t) = at^2 + bt + c$.

In GeoGebra input Curve[t, 0, a*t² +b*t + c, t, ParameterStart, ParameterEnd] where ParameterStart and ParameterEnd give the desired range of values for t.

From (4), let $x(t) = -\frac{b}{2a}$, $y(t) = t$, and $z(t) = a\left(-\frac{b}{2a}\right)^2 - at^2 + b\left(-\frac{b}{2a}\right) + c$.

In GeoGebra input Curve[-b/(2*a), t, a*(-b/(2*a))² - a*t² + b*(- b/(2*a))+c, t, ParameterStart, ParameterEnd] where ParameterStart and ParameterEnd give the desired range of values for t.
This will generate a 3D graph similar to *figure 4*.

### Higher Degree Polynomials

Note when extending this idea to higher degree polynomials, be aware that some choices of coefficients will lead to undefined quantities. These need to be addressed on a case by case basis. For example, the restricted domain of the general quartic is
(9) $y = 0$
or

(11) $y = \pm\sqrt{\frac{4ax^3 + 3bx^2 + 2cx + d}{4ax + b}}$.

The case when the quantity $4ax + b = 0$ must be considered separately.